\documentclass[graybox]{svmult}

\usepackage{type1cm}        % activate if the above 3 fonts are
\usepackage{makeidx}         % allows index generation
\usepackage{graphicx}        % standard LaTeX graphics tool
\usepackage{multicol}        % used for the two-column index
\usepackage[bottom]{footmisc}% places footnotes at page bottom

\usepackage{newtxtext}       % 
\usepackage{newtxmath}       % selects Times Roman as basic font
\usepackage[hidelinks]{hyperref}
\newcommand{\J}{\mathbf{J}}

\usepackage{color} 
\usepackage{psfrag}

\newcommand{\halb}{{\frac{1}{2}}}

\makeindex             % used for the subject index
\begin{document}

\title*{On numerical methods for hyperbolic PDE with curl involutions}
\author{M. Dumbser, S. Chiocchetti and I. Peshkov}
\authorrunning{Dumbser et al.} 
\institute{Laboratory of Applied Mathematics, University of Trento, Via Mesiano 77, 38123 Trento, Italy, \email{michael.dumbser@unitn.it, simone.chiocchetti@unitn.it, ilya.peshkov@unitn.it}
}
\maketitle

\abstract{
	In this paper we present three different numerical approaches to account for curl-type involution constraints in hyperbolic partial differential equations for continuum physics. All approaches have a direct analogy to existing and well-known divergence-preserving schemes 
	for the Maxwell and MHD equations. The first method consists in a generalization of the Godunov-Powell terms, which means adding suitable multiples of the involution constraints to the PDE system in order to achieve the symmetric Godunov form. The second method is an extension of the generalized Lagrangian multiplier (GLM) approach of Munz et al., where the numerical errors in the involution constraint are propagated away via an augmented
	PDE system. The last method is an exactly involution preserving discretization, similar to the exactly divergence-free schemes for the Maxwell and MHD equations, making use of appropriately staggered meshes. We present some numerical results that allow to compare all three approaches with each other. 
}

\section{Introduction} 

Very recently, several novel hyperbolic PDE systems were proposed for the description of dynamic processes in continuum physics that are endowed with curl-type involutions, i.e. where the curl of a certain set of variables either has to vanish or has to assume a prescribed value. The most prominent examples are the system of nonlinear hyperelasticity of Godunov, Peshkov and Romenski (GPR model) \cite{GodunovRomenski72,PeshRom2014,GPRmodel,Rom1998} written in terms of the distortion field $\mathbf{A}$, the conservative compressible multi-phase flow model of Romenski et al. \cite{Rom1998,RomenskiTwoPhase2010}, the new hyperbolic model for surface tension and the  recent hyperbolic reformulation of the Schr\"odinger  equation of Gavrilyuk and Favrie et al.  \cite{Schmidmayer2016,Dhaouadi2018}, as well as first order reductions of the Einstein field equations, such as those proposed, e.g., in \cite{Alic:2009,Brown2012,ADERCCZ4,GLMFOCCZ4}.
Many, but not all, of the aforementioned mathematical models fall into the larger class of symmetric hyperbolic and thermodynamically compatible (SHTC) systems, studied by Godunov and  Romenski et al. in \cite{Godunov1961,Rom1998,Godunov:2003a,SHTC-GENERIC-CMAT}. 
Involution constraints in general are \textit{stationary} differential equations that are satisfied by the governing PDE system \textit{for all times} if they are satisfied by the initial data. 
The most famous involution is the divergence-free condition of the magnetic field in the Maxwell and magnetohydrodynamics (MHD) equations. As a consequence, a lot of research has been dedicated in the past to the appropriate numerical discretization of PDE with divergence constraints. However, much less is known on curl-preserving numerical schemes for PDE with curl involutions. 
In the context of the Maxwell and MHD equations, the most common involution preserving numerical schemes fall into the following three categories: 
\begin{enumerate} 
	\item \textit{Exactly divergence-free schemes}, such as those proposed in \cite{Yee66,BalsaraSpicer1999,BalsaraAMR,GardinerStone,DeVore,balsarahlle2d,balsarahlle3d,ADERdivB}, which make use of the definition of the electromagnetic quantities on appropriately \textit{staggered grids}. To the best knowledge of the authors, the only extensions to curl-type involutions are those presented in \cite{HymanShashkov1997,JeltschTorrilhon2006,Torrilhon2004} so far.  
	\item The formally nonconservative \textit{Godunov-Powell terms}, which go back to a numerical implementation by Powell \cite{PowellMHD1} of the symmetrizing terms of the MHD equations found by Godunov in \cite{God1972MHD}, and which consist in adding suitable \textit{multiples} of the divergence-free condition to the induction, momentum and energy equations. Note that at the analytical level, all these terms are exactly zero, but they are in general non-zero for certain numerical discretizations of the equations that are not in the class of exactly divergence-free schemes. These nonconservative terms which formally correspond to zero were nevertheless needed in order to \textit{symmetrize} the MHD system and to make it at the same time thermodynamically compatible, i.e. to give it the aforementioned SHTC structure. The \textit{obvious disadvantage} of this approach is that it only works for the MHD equations, where a velocity vector is available, since all Godunov-Powell terms are proportional to the velocity. Hence, for the vacuum Maxwell equations, where such a velocity vector does not exist, the approach is not suitable.  
	\item The generalized Lagrangian multiplier (GLM) approach forwarded by Munz et al. in \cite{MunzCleaning,Dedneretal} for the Maxwell and MHD equations. The main idea here consists in solving an \textit{augmented evolution system}, where an artificial scalar cleaning variable $\varphi$ is added and coupled to the induction equation, so that divergence errors in the magnetic field cannot accumulate, but rather propagate away via acoustic-type waves. The advantage of this approach is that it works for MHD as well as for the Maxwell equations and that it does not add any nonconservative terms to the governing equations.  
\end{enumerate} 

At this point, we recall the hyperbolic GLM approach of Munz et al. \cite{MunzCleaning,Dedneretal} in more detail. In this paper we make use of the Einstein summation convention, which implies summation over two repeated indices. We furthermore use the abbreviations  
$\partial_t = \partial / \partial t$,  $\partial_k = \partial / \partial x_k$. The fully anti-symmetric Levi-Civita symbol is denoted by $\varepsilon_{ijk}$. The induction equation in electrodynamics is well-known and reads 
\begin{equation}
\partial_t B_k + \varepsilon_{kij} \partial_i E_j = 0.   
\label{eqn.induction} 
\end{equation} 
Here, $B_k$ and $E_j$ denote the magnetic and the electric field, respectively. An immediate consequence of the induction equation is the 
involution constraint 
\begin{equation}
\mathcal{I} = \partial_m B_m = 0, 
\label{eqn.divBzero}
\end{equation} 
which states that the magnetic field will remain divergence-free for all times, if it was initially divergence-free. As already mentioned above, a classical way to preserve a divergence-free
magnetic field within a numerical scheme is the use of an exactly divergence-free discretization on appropriately  staggered meshes, see
e.g. \cite{Yee66,DeVore,BalsaraSpicer1999,Balsara2004,GardinerStone,ADERdivB}. 
The very popular GLM method proposed by Munz et al. in \cite{MunzCleaning,Dedneretal} is an alternative to exactly constraint--preserving schemes and requires only small changes at the PDE level. Instead of the original induction equation, the following \textit{augmented induction equation} is solved: 
\begin{eqnarray}
\label{eqn.induction.glm} 
\partial_t B_k + \varepsilon_{kij} \partial_i E_j + \textcolor{red}{ \partial_k \varphi } &=& 0, \\   
\label{eqn.phi.glm} 
\textcolor{red}{\partial_t \varphi + a_d^2  \, \partial_m B_m }  &=& \textcolor{red}{-\epsilon_d \varphi}.   
\end{eqnarray} 
Here, $\varphi$ is the new cleaning scalar, $a_d$ is an \textit{artificial} cleaning speed and $\epsilon_d$ is a small damping 
parameter. For convenience, the new terms in the augmented PDE system \eqref{eqn.induction.glm} and \eqref{eqn.phi.glm} with respect 
to the original induction equation \eqref{eqn.induction} are highlighted in red. It is easy to see that for $a_d \to \infty$ 
the equation \eqref{eqn.phi.glm} leads to $\partial_m B_m \to 0$, i.e. in the asymptotic limit the involution constraint \eqref{eqn.divBzero} will be preserved. 
 
In the remaining part of this paper, we will show the natural extensions of the exactly divergence 
free schemes, the Godunov-Powell terms and the GLM cleaning to curl-type involutions. We will show 
computational results for the new hyperbolic surface tension model \cite{Schmidmayer2016} and close 
with some concluding remarks and an outlook to future work.

\section{Model problem and different approaches to account for the curl involution}  
\label{sec.model} 

We illustrate the basic ideas on the following simple toy model, in order to ease notation and to facilitate the understanding of the underlying concepts. Consider the following evolution system for one scalar $\rho$ and two vector fields 
$v_k$ and $J_k$:    
\begin{eqnarray}
\label{eqn.rho}
\partial_t \rho + \partial_i \left( \rho v_i  \right) &=& 0,  \\ 
\label{eqn.vk}
\partial_t (\rho v_k) + \partial_i \left( \rho v_i v_k + \rho^2 E_\rho \, \delta_{ik} + \rho  J_k E_{J_i} \right) &=& 0,  \\ 
\label{eqn.toy} 
\partial_t J_k + \partial_k ( v_m J_m ) + v_m \left( \partial_m J_k - \partial_k J_m \right)  &=& 0. 
\end{eqnarray} 
Here, $E=E(\rho,v_k,J_k)$ is a specific total energy potential and $E_\rho$ and $E_{J_k}$ are the derivatives of the energy potential with respect to the state variables $\rho$ and $J_k$. In particular, $p = \rho^2 E_\rho$ is the fluid pressure. The above system satisfies the additional energy conservation law 
\begin{equation}
 \partial_t (\rho E) + \partial_k \left( v_k (\rho E) + v_i \left(\rho^2 E_\rho \delta_{ik} + \rho J_i E_{J_k} \right)  \right) = 0. 
 \label{eqn.energy} 
\end{equation}
It is easy to see that the PDE \eqref{eqn.toy} 
is endowed with the linear involution constraint  
$\mathcal{I}_{mk} = \partial_m J_k - \partial_k J_m = 0$, i.e. if the curl of $J_k$ is zero for the initial data, then it will remain zero for all times. 
For a general purpose numerical method applied to \eqref{eqn.toy}, it is very hard to guarantee $\mathcal{I}_{mk}=0$ at 
the discrete level. We stress that for smooth solutions, at the continuous level all the following reformulations of the PDE 
system are completely equivalent. The main differences arise at the discrete level. 

\subsection{SHTC structure and Godunov-Powell terms for curl involutions} 

The system \eqref{eqn.toy} can be written in symmetric hyperbolic form \cite{Godunov1961,Rom1998} by adding the term 
$\rho E_{J_k} \left( \partial_i J_k - \partial_k J_i \right) = 0$  
to the momentum equation. Note that the term $v_m \left( \partial_m J_k - \partial_k J_m \right)$ proportional to the
velocity field and to the curl of $J_k$ is already contained in the evolution equation for $J_k$ in order to make the
system Galilean invariant. The modified system then reads 
\begin{eqnarray}
\label{eqn.toy.shtc} 
\partial_t \rho + \partial_i \left( \rho v_i  \right) &=& 0,  \\ 
\partial_t (\rho v_k) + \partial_i \left( \rho v_i v_k + \rho^2 E_\rho \, \delta_{ik} + \rho  J_k E_{J_i} \right)  + \textcolor{red}{\rho E_{J_k} \left( \partial_i J_k - \partial_k J_i \right)} &=& 0,  \\ 
\label{eqn.shtc.Jk} 
\partial_t J_k + \partial_k ( v_m J_m ) + v_m \left( \partial_m J_k - \partial_k J_m \right)  &=& 0, 
\end{eqnarray}
where we have highlighted the additional symmetrizing term in red. 
Introducing the notation $\mathcal{E} = \rho E$, $m_i = \rho v_i$, $r = \mathcal{E}_\rho$, $v_i = \mathcal{E}_{m_i}$, $\eta_i =
\mathcal{E}_{J_i}$, and the Legendre transform $L(\mathbf{p})$ of the potential $\mathcal{E}(\mathbf{q})$ as
\begin{equation}
 L(\mathbf{p}) = \mathbf{q} \cdot \mathcal{E}_{\mathbf{q}} - \mathcal{E} = \rho \mathcal{E}_\rho + m_i \mathcal{E}_{m_i} + J_i \mathcal{E}_{J_i} - \mathcal{E}(\mathbf{q}),   
\end{equation}
with the vector of conservative variables $\mathbf{q} = L_\mathbf{p} = \left( L_r, L_{v_i}, L_{\eta_i} \right) = \left( \rho, m_i, J_i \right)$ and the vector of thermodynamic dual variables $\mathbf{p} = \mathcal{E}_\mathbf{q} = \left( \mathcal{E}_\rho, \mathcal{E}_{m_i}, \mathcal{E}_{J_i} \right) = \left( r, v_i, \eta_i \right)$, one can write the above system \eqref{eqn.toy.shtc}-\eqref{eqn.shtc.Jk} in the \textit{symmetric Godunov form}  
\begin{eqnarray}
    \partial_t L_{r} + \partial_k \left[ \left(v_k L \right)_r \right] &=& 0, \\ 
    \partial_t L_{v_i} + \partial_k \left[\left(v_k L \right)_{v_i}\right] + L_{\eta_i} \partial_k \eta_k - L_{\eta_k} \partial_i \eta_k &=& 0, \\ 
    \partial_t L_{\eta_i} + \partial_k \left[ \left(v_k L \right)_{\eta_i} \right] - L_{\eta_i} \partial_k v_k + L_{\eta_k} \partial_i v_k &=&0. 
\end{eqnarray}   
The modified system \eqref{eqn.toy.shtc}-\eqref{eqn.shtc.Jk} is not only symmetric hyperbolic for convex potentials $L$, but is also numerically much better behaved concerning the curl involution on $J_k$ when solved with a general purpose scheme. 

\subsection{GLM curl cleaning}  

As already mentioned before, the main advantage of the GLM approach of Munz et al. 
\cite{MunzCleaning,Dedneretal} for divergence constraints is its ease of implementation and the fact that it does not necessarily require a velocity field, since the transport of the divergence errors is achieved via acoustic-type waves. Here, in the case of curl involutions, we add a Maxwell-type subsystem, i.e. curl errors propagate away via electro-magnetic-type waves. The disadvantage of GLM curl cleaning is the need to add a rather large number of auxiliary evolution quantities to the system. 
The GLM curl cleaning proposed in \cite{GLMFOCCZ4,SHTCSurfaceTension} can  be explained on the toy system  \eqref{eqn.rho}-\eqref{eqn.toy} as follows. The original governing PDE system \eqref{eqn.rho} - \eqref{eqn.toy} is 
simply \textit{replaced} by the following \textit{augmented system} that accounts for the curl constraint on $J_k$: 
\begin{eqnarray}
\label{eqn.toy.glm} 
\partial_t \rho + \partial_i \left( \rho v_i  \right) &=& 0,  \\ 
\partial_t (\rho v_k) + \partial_i \left( \rho v_i v_k + \rho c_0^2 J_i J_k \right) &=& 0,  \\ 
\label{eqn.glm.Jk} 
\partial_t J_k + \partial_k ( v_m J_m ) + v_m \left( \partial_m J_k - \partial_k J_m \right) + \textcolor{blue}{\varepsilon_{klm} \partial_l \psi_m} &=& 0, \\ 
\label{eqn.glm.psi} 
\textcolor{blue}{ 
	\partial_t \psi_k - a_c^2 \, \varepsilon_{klm} \partial_l J_m } \textcolor{red}{ + \partial_k \varphi } &=& 
\textcolor{blue}{  - \epsilon_{c} \, \psi_k,  }  \\  
\label{eqn.glm.phi} 
\textcolor{red}{ 
	\partial_t \varphi + a_d^2 \, \partial_m \psi_m } &=& \textcolor{red}{- \epsilon_{d} \varphi}, 
\end{eqnarray}
where $a_c$ is a new cleaning speed associated with the curl cleaning. The new terms associated with the curl cleaning are highlighted in blue, for convenience, while the terms of the original PDE \eqref{eqn.toy} are written in black. Since the evolution equation for the cleaning vector field $\psi_k$ has formally the same structure as the induction equation \eqref{eqn.induction} of the Maxwell equations, it is again endowed with the divergence-free constraint 
$\partial_m \psi_m = 0$, which is taken into account via the classical GLM method (red terms).  
It is easy to see that from \eqref{eqn.glm.psi} for 
$a_c \to \infty$ we obtain $ \epsilon_{klm} \partial_l J_m \to 0$ in the limit, thus satisfying 
the involution in the sense $\mathcal{I}_{mk} \to 0$. The augmented system \eqref{eqn.toy.glm}-\eqref{eqn.glm.phi} can now be solved with any standard numerical method for nonlinear systems of hyperbolic partial differential equations. The main advantage over the Godunov-Powell terms proposed in the previous section is the fact that the GLM curl cleaning does \textit{not} destroy conservation of  momentum and it also works in absence of a physical velocity field $v_k$.

\subsection{An exactly curl-free discretization} 

\begin{figure}[!htbp]
	\begin{center} 
	\begin{tabular}{cc}
		\includegraphics[width=0.452\textwidth]{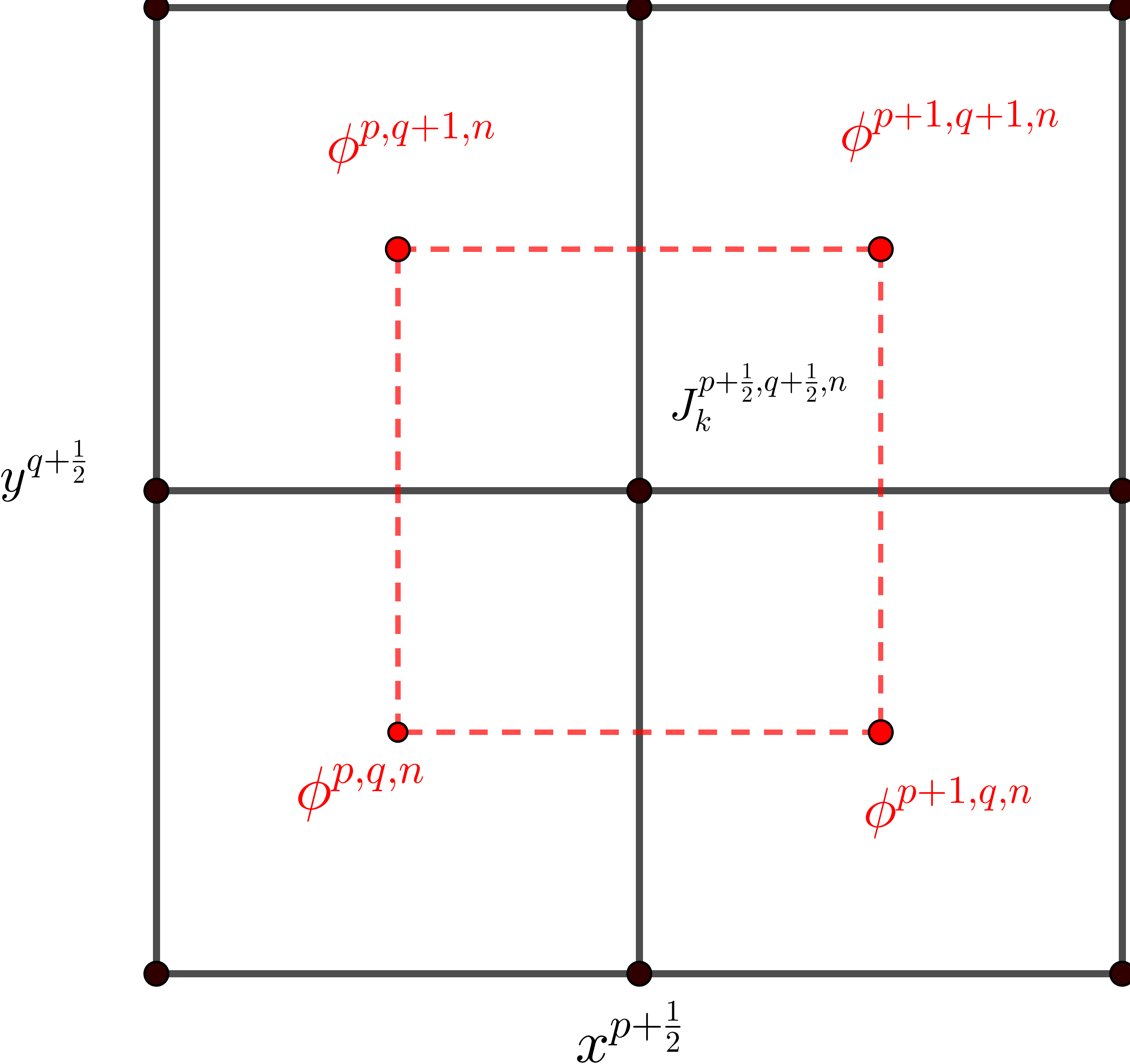} & 
		\includegraphics[width=0.44\textwidth]{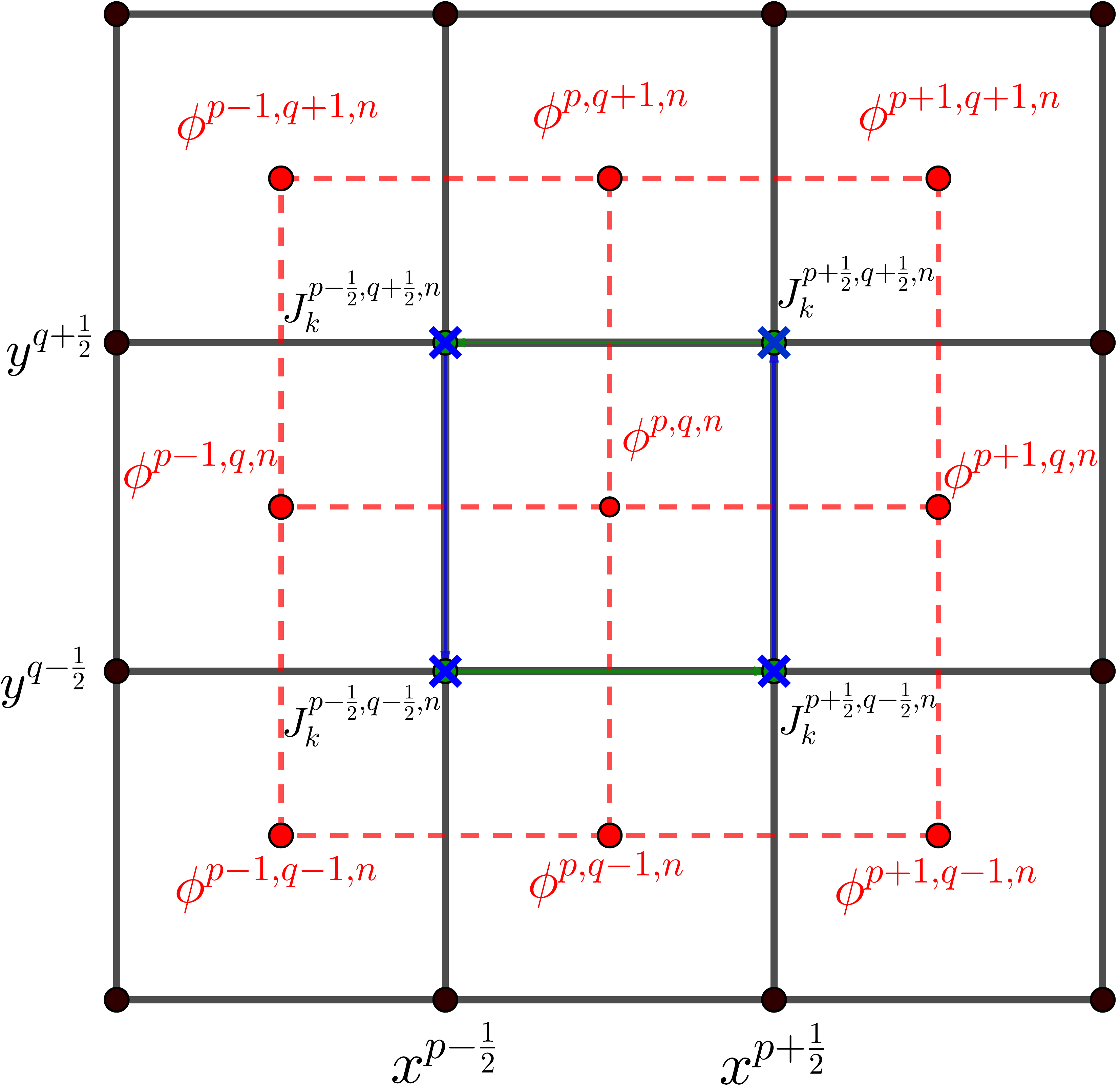} 
	\end{tabular}    
	\end{center} 	
	\caption{Left: stencil of the discrete gradient operator, which computes the corner gradient of a scalar field defined in the  barycenters. Right: stencil of the discrete curl operator, defining a curl inside the barycenter using the vector field in the corners. }  
	\label{fig.grad.curl}
\end{figure}

Here we present a compatible discretization that satisfies the curl constraint \textit{exactly} at the discrete level. For this 
purpose, we use an appropriately \textit{staggered mesh}, with the field $J_k$ defined in the vertices of the main grid and the 
scalar field $\phi = v_m J_m$ defined in the barycenters of the primary control volumes. 
To avoid confusion between tensor indices and discretization indices, throughout this paper we will use the \textit{subscripts}  $i,j,k,l,m$ for tensor indices and the superscripts $n,p,q,r,s$ for the \textit{discretization indices} in time 
and space, respectively. 
The discrete spatial coordinates will be denoted by $x^p$ and $y^q$, while the set of discrete times will be denoted by $t^n$.  
The $z$ component of the \textit{discrete curl} $\nabla^h \times $ of a discrete vector field 
$\J^{h,n}$
is denoted by $\left( \nabla^h \times \J^{h,n} \right) \cdot \mathbf{e}_z$ and its degrees of freedom are naturally defined as
\begin{eqnarray}
\left( \nabla^{p,q} \times \mathbf{J}^{h,n} \right) \cdot \mathbf{e}_z &=&  
\halb \frac{J_2^{p + \halb,q +\halb,n} + J_2^{p + \halb,q -\halb,n} - J_2^{p - \halb,q + \halb,n} - J_2^{p - \halb, q -\halb,n} }{\Delta x} - 
\nonumber \\ &&  
\halb \frac{J_1^{p + \halb,q +\halb,n} + J_1^{p - \halb,q +\halb,n} - J_1^{p + \halb,q - \halb,n} - J_2^{p - \halb, q -\halb,n} }{\Delta y} 
\label{eqn.rot} 
\end{eqnarray}
making use of the \textit{vertex-based} staggered values of the field $\J^{h,n}$, see the right panel in Fig. \ref{fig.grad.curl}. In Eqn. \eqref{eqn.rot} the symbol $\epsilon_{ijk}$ is the usual Levi-Civita tensor. Eqn. \eqref{eqn.rot} defines a discrete
curl on the control volume $\Omega^{p,q}$ via a discrete form of the Stokes theorem based on the trapezoidal 
rule for the computation of the integrals along each edge of $\Omega^{p,q}$. 
Last but not least, we need to define a discrete gradient operator that is compatible with the discrete curl,
so that the continuous identity
\begin{equation}
\nabla \times \nabla \phi = 0
\label{eqn.rotgrad} 
\end{equation}
also holds on the discrete level. If we define a scalar field in the barycenters of the control volumes $\Omega^{p,q}$ as
$\phi^{p,q,n}=\phi(x^p,y^q,t^n)$ then the corner gradient generates a natural discrete gradient operator $\nabla^{h}$ 
of the discrete scalar field $\phi^{h,n}$  that defines a discrete gradient in all vertices of the mesh. 
The corresponding degrees of freedom generated by $\nabla^{h} \phi^{h,n}$ read 
\begin{equation}
\nabla^{p+\halb,q+\halb}  \phi^{h,n} = \partial_k^{p+\halb,q+\halb} \phi^{h,n} = \left( \begin{array}{c}  
\halb \frac{\phi^{p + 1, q + 1,n} + \phi^{p + 1, q,n} - \phi^{p, q + 1,n} - \phi^{p, q,n} }{\Delta x}  \\
\halb \frac{\phi^{p + 1, q + 1,n} + \phi^{p, q + 1,n} - \phi^{p + 1, q,n} - \phi^{p,q,n} }{\Delta y}   \\ 
\displaystyle 0
\end{array} \right), 
\label{eqn.grad} 
\end{equation}
see the left panel of Fig. \ref{fig.grad.curl}. 
It is then straightforward to verify that an immediate consequence of \eqref{eqn.rot} and \eqref{eqn.grad} is 
\begin{equation}
\nabla^h \times \nabla^h \phi^{h,n} = 0, 
\label{eqn.discrete.curl} 
\end{equation}
i.e. one obtains a discrete analogue of \eqref{eqn.rotgrad}. With this compatible discretization, Eqn. \eqref{eqn.toy}, which contains 
a gradient and a curl operator, can be discretized so that $J_k$ remains curl-free for all times.

\subsection{Numerical results} 

\begin{figure}[!htbp]
	\includegraphics[width=1.0\textwidth]{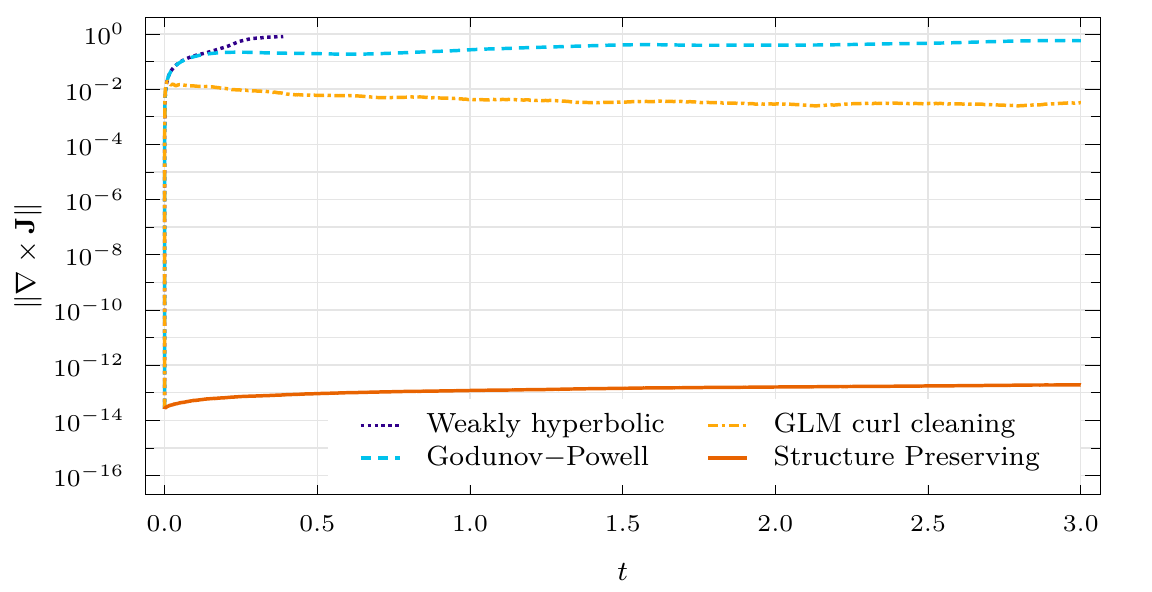} 
	\caption{Temporal evolution of the curl errors for different numerical methods applied to the hyperbolic surface tension model proposed in \cite{Schmidmayer2016}. }  
	\label{fig.curlerror}
\end{figure}

Here we present some numerical results obtained with the three approaches mentioned above, applied to the hyperbolic surface tension model 
of Gavrilyuk et al. \cite{Schmidmayer2016}. It is important to note that the original model \cite{Schmidmayer2016} is only \textit{weakly hyperbolic} and thus \textit{not} suitable for a stable numerical discretization with a general purpose scheme. In Fig. \ref{fig.curlerror} we compare the numerical results obtained for the original weakly hyperbolic model, for the non-conservative Godunov-Powell terms, for the GLM curl cleaning and for the exactly curl-free discretization. The results clearly show that the weakly hyperbolic system becomes unstable with a general purpose scheme, while the non-conservative Godunov-Powell terms allow a stable discretization. Even better results are obtained for the conservative GLM curl cleaning approach. The best results are obtained by the exactly curl-free (structure preserving) scheme, which can be even directly applied to the weakly hyperbolic system, thus emphasizing the important role of the curl involution at the continuous and discrete level. 

\section{Conclusions}
\label{sec.concl} 
We have outlined three possible extensions of divergence-free schemes to hyperbolic PDE systems with curl-type involution constraints, namely i) the classical Godunov-Powell approach based on the symmetrization of the governing PDE system, ii) the hyperbolic GLM cleaning approach that accounts for the involution constraint via an augmented PDE system and in which the numerical errors of the involution are transported away via a Maxwell-type subsystem and iii) an exactly curl-free scheme based on appropriately staggered meshes.  
Future work will consist in an extension of the exactly curl-free approach to higher order of accuracy and the application to 
other PDE systems with curl-type involutions, such as those presented in \cite{Schmidmayer2016,RomenskiTwoPhase2007,RomenskiTwoPhase2010,Schmidmayer2016,Dhaouadi2018}. 
First preliminary results of the authors indicate that the use of exactly curl-free schemes for hyperbolic PDE systems with curl involutions are by far superior in performance and accuracy compared to the Godunov-Powell terms and compared to the GLM cleaning approach.

\begin{acknowledgement}
The research presented in this paper has been funded by the European Union's Horizon 2020 Research and Innovation Programme under the project \textit{ExaHyPE}, grant no. 671698 and by the Deutsche Forschungsgemeinschaft (DFG) under the project DROPIT (Droplet Interaction Technologies), grant no. GRK 2160/1. 
MD also acknowledges financial support from the Italian Ministry of Education, University and Research (MIUR) via the Departments of Excellence Initiative 2018--2022 attributed to DICAM of the University of Trento (grant L. 232/2016) 
and via the PRIN 2017 project.  

This is a pre-print of the following work: 
{G.V. Demidenko, E. Romenski, E.F. Toro, M. Dumbser (Eds.)}, 
\textit{``Continuum Mechanics, Applied Mathematics and Scientific Computing: Godunov's Legacy''}, 2020, 
Springer International Publishing. 
Reproduced with permission of Springer Nature Switzerland AG.
\texttt{DOI: 10.1007/978-3-030-38870-6}.
\end{acknowledgement}

\bibliographystyle{plain}
\bibliography{./GLMRot}

\end{document}